\documentclass{article}

\include{8bitdefs}
\usepackage[applemac]{inputenc}
\usepackage[french]{babel}

\title{Le Théorème de Riesz-Raikov-Bourgain pour un endomorphisme
algébrique de $\R^p$ }

\author{{par Jean-Claude LOOTGIETER}\\
\\
{ Laboratoire de Probabilités et Modèles Aléatoires}\\
{ Université Pierre et Marie Curie}\\
{ 4, Place Jussieu}\\
{ 75252 PARIS Cedex 05, France} }

\usepackage{amsmath,amsfonts}
\usepackage{amssymb}
\begin{document}

\newcommand{\ind}{\mbox{\rm 1\hspace{-0.032in}I}}
\newcommand{\valabs}[1]{\lvert #1 \rvert} %$\valabs{x}
\newcommand{\norm}[1]{\lVert #1 \rVert} %$\norm{x}$
\newcommand{\N}{\mathbb{N}}
\newcommand{\Z}{\mathbb{Z}}
\newcommand{\R}{\mathbb{R}}
\newcommand{\C}{\mathbb{C}}
\newcommand{\Q}{\mathbb{Q}}
\newcommand{\T}{\mathbb{T}}
\newcommand{\eps}{\varepsilon}

\newtheorem{lem}[equation]{Lemme}%[subsection]
\newtheorem{rem}[equation]{Remarque}%[subsection]
\newtheorem{theo}[equation]{Théorème}%[section]
\newtheorem{cor}[equation]{Corollaire}%[section]
\newtheorem{guess}[equation]{Sous-lemme}%\newtheorem{note}{Note}[subsection]

\hspace*{-0.3in}

\noindent   {\bf {SUR  L'\' {E}NONC\' {E}  D'UN TH\' {E}OR\`{E}ME CONCERNANT LES\\
OP\' {E}RATEURS  POSITIFS  SUR LES ESPACES  $\bf{\it{L}}_{ \it{{p}}}\,\,(1< \it{p}<\infty)$  \\  DONT LA SUITE DES
PUISSANCES EST SOUS-ADDITIVE}}

%\noindent   {\bf {Sur l' énoncé incorrect d'un théorème concernant les opérateurs positifs  sur les espaces $\bf{L^p\,\,(1<p<\infty)}$  dont  la suite des puissances est sous-additive.}}\\
\vspace{1 cm}
\hspace{3 cm}{J.\,\,C. LOOTGIETER}\\
%\noindent Note de Lootgieter Jean-Claude. Université Pierre et Marie Curie.\\

 \vspace{1 cm}
 \noindent{ABSTRACT}$-$In this article we give a counter-example on  the statement of a theorem appearing  in a note  $[3]$ of  A. Brunel concerning the study of  positive operators on the the spaces   $L_p\,\,(1<p<\infty)$  which the sequence of the powers is sub-additive.\\
  \vspace{0,5 cm}
 
\noindent $Introduction.$ \,\,\,Dans un article paru en 1992 $[2] $ A. Brunel énonce un théorème généralisant   le  théorème suivant démontré par  A. Akcoglu $[1]$: soit 
\begin{equation}
M_n(T)=\frac{1}{n}\sum_{j=0}^{n-1} T^j
\end{equation} 
les moyennes de Césaro des puissances d'une contraction linéaire positive sur l'espace $L_p\,\,(1<p<\infty)$
d'une mesure $\sigma$-finie; alors, pour toute $f \in$ $L_p$, on a l'inégalité (dite ergodique dominée)
\begin{equation}
\norm {\sup_{n}\vert  M_n(T)(f)\vert\,}_p\leq \frac{p}{p-1} \,\norm f_{p}.
\end{equation} 
Cette inégalité permet d'en déduire la convergence p.p (p.p=presque-partout) des moyennes de Césaro
$M_n(T)(f)$ pour $f\in L_p$.\\

\noindent Dans $[2] $ on suppose   que  $T$ n'est pas nécessairement une  contraction, mais que  les puissances  
 $T^n$ de $T$ sont à moyennes bornées:
\begin{equation}
\sup_{n}\,\norm {M_n(T)}_{p}< \infty
\end{equation} 
et  on démontre (c.f $[2]$, théorème 1)  l'inégalité ergodique dominée: il existe une constante $c$  finie (a priori dépendant de $T$) telle que, pour toute $f \in$ $L_p$, 
\begin{equation}
 \norm {\sup_{n}\vert  M_n(T)(f)\vert}_p\leq c  \,\norm f_{p}.
 \end{equation}
  
 \noindent Dans la suite $c$, $C$, $c_0$, $C_0$, etc... désigneront des constantes $>0$  et finies variées.\\
 
\noindent Dans $[2]$ la démonstration de (4) utllise  le barycentre des puissances $T^j$ de $T$:
\begin{equation}
A:=\sum_{j=0}^{\infty} \alpha_{j}T^j
\end{equation} 
où $\alpha$ est la probabilité sur $\N$ de fonction génératrice 
\begin{equation}
f(x):=\frac{1-\sqrt {1-x}}x.
\end{equation} 

\noindent De l'expression des convoluées $\alpha^{n}$, $n\in \N^*$, de la probabilité $\alpha$:
\begin{equation}
\alpha^n_j=\frac{n}{2(j+n)}2^{1-2j-n}\begin{pmatrix} 2j+n-1\\ j \end{pmatrix}
\end{equation}
découle (c.f $[4]$) que  $A$ est bien un opérateur positif  de  $L_p$ dans  $L_p$ et que ses puissances 
\begin{equation}
A^n=\sum_{j=0}^{\infty} \alpha^n_{j}T^j
\end{equation}
vérifient
\begin{eqnarray}
%\,\,\,&\,\,\,&  \norm{A_n}_p/n \,\,\mathrm{tend \,\,vers }\,\,0.\nonumber \\
(a)&\,\,\,&\sup_n\,\norm{A^n}_p\leq C \sup_n \,\norm{M_n(T)}_p\,, \nonumber  \\
(b)&\,\,\,&\mathrm{pour\,\,tout \,\,n} \in \N,\,\,M_n(T)\leq c\, M_{q(n)}(A),  \nonumber   \\
&\,\,\,&
\mathrm{si}\,\,q(n)=\mathrm{Ent} (\sqrt n)+1,
 \nonumber \\
(c)&\,\,\,& \mathrm{pour \,\,tous\,\, n,\,m\,}\in \N, \,A^{n+m}\leq A^{n}+A^{m}.
\end{eqnarray}

\noindent La propriété de sous-additivité (9)-(c) résulte du fait que, pour tout $j$ fixé, la suite $\alpha^n_{j}/n$, $n\in\N^*$, est décroissante 
et donc que 
$$\mathrm{pour}\,\,\,\mathrm{tous }\,\,n,\,m \in\N,\,\,\alpha^{n+m}_{j} \leq \alpha^n_{j}+\alpha^m_{j}.$$\\
\noindent  De (9)-(a) et (3) résulte que les puissances %$A^n$
de $A$ sont bornées en norme:
\begin{equation}
\sup _n \,\norm {A^n}_p< \infty.
\end{equation}
\noindent $Remarque\, \it{1}$- De (10) résulte que $\lim_n \norm{A_n}_p/n=0$ et donc, 
pour toute $f\in L_p$,  la suite $A^n(f)/n$ converge dans $L_p$ vers $0$;  d'autre part, pour $f\in L_p$,
 $f\geq 0$, la suite $A^n(f)/n$, $n\in N^*$, est positive, décroissante (puisque, pour tout $j$ fixé, la suite $\alpha^n_{j}/n$, $n\in\N^*$ est décroissante) et  donc converge p.p. vers $0$. \\
 
 \noindent $Remarque \,\it{2}$- Considérons un opérateur positif $A$ de  $L_p$ dans  $L_p$ dont la suite des puissances $A^n$ vérifie la $seule$ propriété de sous-additivité: 
\begin{equation}
 \hspace{2,2 cm}(c)\,\,\,\,\,\,\,\,\,\,\,\,\mathrm{pour \,\,tous\, \,n,\,m\,}\in \N, \,A^{n+m}\leq A^{n}+A^{m},   \hspace{2,2 cm} (9) \nonumber\\
\end{equation}
  La condition (9)-(c)  implique que, $\mathrm{ pour \,\,toute}$ $f \in L_p$, $f\geq 0$, la suite $A^n(f)$ est positive, sous-additive et donc la suite $A^n(f)/n$ converge p.p et, étant dominée par la fonction $A(f)$, converge également dans $L_p$. Pour que cette limite soit nulle il suffit que $\lim_n \norm{A_n}_p/n=0$, \\

\noindent Dans $[2]$, p. 197, une étape essentielle  de la démonstration de l'inégalité (4) est  l'inégalité ergodique dominée:
pour toute $f \in$ $L_p$,  
\begin{equation}
\norm{\sup_n\,\lvert A^n(f)\lvert}_p\leq C\, \,\norm f_{p},
\end{equation}
condition qui implique que les  $ A^n$ sont bornées en norme, i.e. (10).\\

\noindent Ryotaro Sato $[7]$ (cité dans $[3]$) souligne que les arguments développés dans $[2]$ 
pour aboutir à (11) s'appuient {\it seulement} sur la propriété de sous-additivité (9)-(c) et donne le contre-exemple d'un type d'opérateurs
positifs A vérifiant  (9)-(c), mais dont les normes des puissances de $A$ %$ A¨n$ 
ne sont pas bornées:
\begin{equation}
A=I+B, \mathrm{avec} \,\,\,  B> 0, \,\, \mathrm{et}\,\, B^2=0.
\end{equation}
De (12) résulte que $A^n=I+nB$, donc $A^{n+m} \leq A^n+A^m$ et  $\norm{A^n}_p> n\norm{B}_p$.\\
\noindent Par exemple, dans le cas où $L_p=l^p(\{0,1\})=\R^2$ avec $\mu=(1, 1)$, c'est le cas de l'opérateur $A$ dont la matrice dans la base canonique de $\R^2$ est
\begin{equation}
A=\begin{pmatrix}1&&&a\\  0&&&1\end{pmatrix} \,\,\, \mathrm{et \,\,\,donc}\,\,\,A^n=\begin{pmatrix}1&&&na\\  0&&&1\end{pmatrix}
\end{equation}
avec $a>0$.\\

 \noindent A la suite de ce  contre-exemple, A. Brunel dans $[3]$, p. 207, 
 ajoute à la condition (9)-(c) la condition  supplémentaire: 
\begin{equation}
\lim_n \frac{A^n}{n}=0..
\end{equation}
et énonce le théorème suivant (théorème 2 dans $[3]$): \\ 

\noindent{\large{\it{Théorème} }}- {\it { Soit A un opérateur positif sur}}  $L_p$ $(1<p<\infty)$ {\it {dont la suite des puissances est sous-additive. }}  {\it {Si la suite $A^n/n$ $converge \,\, vers \,\,0$}},
$ alors$\\
 {\it {(i) les puissances $A^n$ de $A$  sont bornées en norme,}}\\
 {\it{(ii) pour toute $f\in L_p$, la suite  $A_n(f)$  converge p.p.}}\\
 
% \noindent Le commentaire qui suit l'énoncé de ce théorème (cf.  $[3]$, p. 207)
% assure que les  résultats de $[2]$ sont alors valides.\\

\noindent Nous donnons dans  ce qui suit  un contre-exemple pour lequel les assertions (i) et (ii) de ce 
théorème ne sont pas réalisées.\\
 
\noindent D'abord une remarque:\\
\noindent $Remarque\, \it{3}$- La condition (14)
  est ambiguê. S'agit-il, pour toute $f\in L_p$, de la convergence  $L_p \,\,\mathrm{et \,\,p.p}$ vers $0$ de la suite  ${A^n(f)}/{n}$
 (cf. [3] p. 206) ou de la condition plus forte  (cf. $ remarque \,{\it2}$)
 $\lim_n\norm{A^n}_p/n=0$? 
  Dans la suite nous garderons l'hypothèse 
 $\lim_n\norm{A^n}_p/n=0$.\\
%\noindent Le commentaire qui suit l'énoncé du théorème 2 (cf. [3] p. 207) assure que les  résultats de $[2]$ sont alors valides.\\

%\noindent Nous donnons dans  ce qui suit  un contre-exemple (avec la condition la plus forte\\
 %$\lim_n\norm{A^n}_p/n=0$ ) pour lequel les assertions (i) et (ii) ne sont pas réalisées.\\

\noindent Nous utiliserons le langage et les notations probabilistes. (P pour probabilité, E pour espérance).\\
$X_1,X_2,\ldots X_n$  désigne une suite de variables aléatoires indépendantes à valeurs $\N$ de même loi $\alpha$ ($\alpha$ est la probabilité donnée par (6)). On pose $S_0=0$ et, pour $n\in \N^*$,
 $S_n=X_1+X_2+\ldots+X_n$.\\
\noindent Rappelons d'abord les propriétés de la probabilité $\alpha$ que nous utiliserons: 
 \begin{eqnarray}
(a) &\,\,\,&\mathrm{pour}\,\,\mathrm{tous }\,\,\,j,n,m \in \N,\,\,\,\alpha^{n+m}_{j} \leq \alpha^n_{j}+\alpha^m_{j}, \nonumber \\
%(b)&\,\,\,\,&  \mathrm{quand}\,\, j \rightarrow \infty,  \alpha_j\sim c\,j^{-3/2},\,\,\mathrm{et\,\,\, donc }\,\,\, \alpha_{j}
%\sim \alpha_{j+1} \nonumber \\
(b) &\,\,\,&\mathrm {pour\,\,toute\,\,constante}\,\, c_0\geq 1\,\,\mathrm{il} \,\,\,\mathrm{existe} \,\,\,\mathrm{une\,\,\, constante}\,\,\, C_0  \,\,\,  \nonumber \\
 \,\,&\,\,\,\,&\mathrm {telle} \,\,\,\mathrm {que} 
  \nonumber \\
 \,\,&\,\,\,&\mathrm {si}\,\, c_0\,j\geq n^2,\,\mathrm {alors}  \,\,\alpha_j^n \geq  C_0\, n \,\alpha_j
\end{eqnarray}
Pour la propriété (15)-(b) se reporter à $[6]$, p. 211.\\

\noindent Considérons l' espace:
\begin{eqnarray}
L_p&:=& l^p(\N,\,\alpha) \nonumber \\
&=&\lbrace f: \N \mapsto \R \mid E\big(\vert f(X_1)\vert^p \big)<\infty \rbrace
\end{eqnarray} 
muni de la norme $\norm{f}_p=\big(E\big(\vert f(X_1)\vert^p \big)\big)^{1/p}$.\\

\noindent Pour  $f\in L_p$ et $k\in \N$ posons 
\begin{equation}
A(f)(k):=E\big(f(X_1+k)\big).
\end{equation} 
%D"après  (b)-(16), pour tout $k$,  $\alpha_j \sim \alpha_{j+k}$ $(j\rightarrow \infty)$ et donc\\
%$f(X_1+k)$ est intégrable. $A(f)(k)$ est donc bien défini.
%Formellement $A$ est le $\alpha$-barycentre des itérées de l'opérateur de translation:
%pour $k\in \N$, 
%$T(f)(k)=f(1+k)$, i.e
%\begin{equation}
%A:=\sum_{j=0}^{\infty} \alpha_{j}T^j.
%\end{equation}
%Mais pour la suite nous prendrons pour point de départ  la définition probabiliste (17) de $A$ .\\
%Comme
On a
 \begin{eqnarray}
\sum_{k}^{\infty}\alpha_k\big (E\big(\vert f(X_1+k)\vert \big)\big)^p&\leq& 
\sum_{k}^{\infty} \alpha_k E\big(\vert f(X_1+k)\vert^p \big) \nonumber \\ 
 &=&E\big(\vert f(X_1+X_2)\vert^p \big) \nonumber \\
 &=& \sum_{k=0}^{\infty} \alpha_k^2\vert f(k)\vert^p \nonumber \\
 &\leq & 2 \sum_{k=0}^{\infty} \alpha_k \vert f(k)\vert^p  \mathrm{\,\,\,\,(cf.\,\,(15)-(a))},
    \nonumber \\
 &=& 2 \norm{f}_p^p.
\end{eqnarray}
D'après (18), les $A(f)(k)$ sont donc bien définis et $A$ définit bien un opérateur positif de $L_p$
dans $L_p$ pour lequel
\begin{eqnarray}
 \norm{A(f)}_p
  &\leq & 2^{1/p}\norm{f}_p.
\end{eqnarray}
Une simple récurrence assure que les puissances de $A$ sont données par
 \begin{eqnarray}
A^n(f)(k)&=&E\big(f(S_n+k)\big) \nonumber \\
\,\,\,\,\,\,\,\,&=& \sum_{j=0}^{\infty} \alpha_j^n f(j+k).
\end{eqnarray}
De (15)-(a) et (20 découle que la suite $A^n$ est sous-additive:\\
\begin{equation}
\hspace{2,2 cm}(c)\,\,\,\,\,\,\,\,\,\,\,\,\mathrm{pour \,\,tous\,\, n,\,m\,}\in \N, \,A^{n+m}\leq A^{n}+A^{m},   \hspace{2,2 cm} (9) \nonumber\\
\end{equation}
\noindent De (20) on  déduit 
 \begin{eqnarray}
\norm{A^n(f)}_p^p&=&\sum_{k}^{\infty} \alpha_k\big\vert E\big(f(S_n+k) \big)\vert^p \nonumber \\
 &\leq&\sum_{k}^{\infty} \alpha_k E\big(\vert f(S_n+k)\vert^p \big)  \nonumber \\ 
 &=&E\big(\vert f(S_n+X_{n+1})\vert^p \big) \nonumber \\
 &=& \sum_{k=0}^{\infty} \alpha_k^{n+1}\vert f(k)\vert^p \nonumber \\
 &\leq & (n+1) \sum_{k=0}^{\infty} \alpha_k \vert f(k)\vert^p  \mathrm{\,\,\,\,(cf.\,\,(15)-(a))}
    \nonumber \\
 &=& (n+1) \norm{ f}_p^p.
\end{eqnarray}
Par conséquent $\norm{A^n}_p\leq  (n+1)^{1/p}$ et par suite  $\lim_n \norm{A_n}_p/n=0$.\\
%Compte tenu de la $remarque${\it 2} on conclut que  $A^n(f)/n$ converge dans $L^2$ et p.p. vers 0.\\

\noindent Les conditions  de l'énoncé du théorème sont donc réalisées.\\

\noindent $Remarque\,{\it 4}$- On peut remarquer que $A$ (cf. (17) et (5)) peut s"interpréter comme le barycentre
  $A=\sum_j \alpha_jT^j$ des itérées de l'opérateur de translation $T$: pour tout $k \in \N$, 
  $T(f)(k)=f(1+k)$. \\

\noindent 1) Comme $P(X_1=0)<1$, il est clair que $P(\lim_nS_n=\infty)=1$. Soit la fonction
sur $\N$ défini par $f(k)=k^{2/5p}$. Comme $\alpha_k \sim c\,k^{-3/2}$ il est clair que $f\in L_p$.
On a, puisque $\lim_kf(k)=\infty$, $P(\lim_n f(S_n)=\infty)=1$. Le lemme de Fatou implique que
 $\lim _nE\big(\big(f(S_n)\big)=\infty$. Comme 
 \begin{eqnarray}
 \norm{A^n(f)}_p^p\geq \alpha_0 \big(E\big(f(S_n)\big)\big)^p
%\norm{A^n(f)}_p^p&\geq& \big(E\big(f(S_n)\big)\big)^p \,\,\, \mathrm {puisque} \,\, f\,\, \mathrm {croissante} 
\end{eqnarray}
$\lim_n  \norm{A^n(f)}_p=\infty$. On en déduit que les $\norm{A^n}_p$ ne sont  pas bornées.\\
\noindent   Il en serait de même pour les moyennes de Césaro $\frac{1}{n}\sum_{j=0}^{n-1}A^j(f)$ avec la même fonction $f$.\\

\noindent On peut estimer inférieurement les  $\norm{A^n}_p$. Posons  
\begin{equation}
f_n(k)=0\,\, \mathrm{si}\,\, k< n^2,\,\, =1\,\,\mathrm{si}\,\, k\geq n^2.\\ 
\end{equation}
De  l'équivalence $\alpha_k \sim c k^{-3/2}$ découle que
\begin{equation}
\norm{f_n}_p^p=\sum_{k=n^2}^{\infty}\alpha_k  \sim c n^{-1}.
\end{equation}
De (15)-(b) découle que
\begin{eqnarray}
  \norm{A_n(f_n)}_p^p)&=& 
  \sum_{k=0}^{\infty}\alpha_k\big(E\big(\big(f_n(S_{n}+k)\big)\big)^p \nonumber \\
 &\geq&\big(E\big(\big(f_n(S_{n})\big)\big)^p \nonumber \\ 
 &=&\bigg(\sum_{j=n^2}^{\infty}\alpha_j^n\bigg)^p \nonumber \\
  &\geq&\bigg (C\, n\sum_{j=n^2}^{\infty}\alpha_j\bigg)^p \nonumber \\
 &\geq& c_1.
 \end{eqnarray}
\noindent  Compte tenu de (24) on en déduit que $ \vert \lvert A^n\vert \lvert_p \geq  C n^{1/p}$. \\

\noindent 2)  Comme $P(\lim_nS_n=\infty)=1$ il est très facile de définir des fonctions $f \in L_p$
telles que la suite $A^n(f)(k)$ ne converge en aucun point $k\in\N$ tout en étant bornée ou non. L'exemple $f(k)=k^{2/5p}$ du 1) précédent en est déjà un exemple: $\lim_n A(f)(k)=\lim E(f(S_n+k))=\infty$. \\
 Il est également très facile de définir des fonctions $f \in L^p$ telles que la suite des moyennes de Césaro $\frac{1}{n}\sum_{j=0}^{n-1}A^j(f)(k)$ ne converge en aucun point $k\in\N$ tout en étant bornée ou non. \\
 
\noindent Les assertions (i) et (ii) du théorème ne sont donc pas réalisées.\\

\noindent $Remarque\, \it{5}$- Pour la preuve proposée dans $[3]$ de ce  théorème l'auteur se borne à indiquer d'appliquer à $A$  la  
méthode développée dans $[2]$.\\
L'étape essentielle de la méthode développée dans $[2]$ et appliquée à $A$ est l'inégalité ergodique
dominée ( p. 197): pour toute $f \in$ $L_p$
$$\norm{\sup_n\lvert A^n(f)\lvert}_p\leq C\, \,\norm f_{p}, \,\,\,\,\,\,(*)$$\\
condition qui implique l'assertion (i), assertion que notre exemple  contredit.\\

\noindent C'est la même inégalité $(*)$, en tenant compte de l'inégalité (9)-(b), qui
conduit à  l'inégalité ergodique dominée (4) pour les moyennes de Césaro de $T$ sans supposer explicitement  que $T$ soit a priori Césaro-borné, mais en supposant seulement  que la définition (5)  de $A$ définit bien un opérateur de $L_p$ dans  $L_p$ (cf. remarque (c) dans $[2]$, p. 207).\\
Notre exemple en fournit un contre-exemple (cf. $ remarque\,\it {4}$) car il est
très facile de remarquer que les moyennes de Césaro de l'opérateur de translation $T$ ne satisfont pas
 l'inégalité ergodique dominée (4).\\
 \,\\
 \,\\

%\noindent  {\it{\bf Conclusion}} \\
%\noindent  Nous concluons %qu'une lecture critique  de $[2]$  est souhaitable 
%et 
%que la  preuve de l'inégalité ergodique dominée (4) sous l'hypothèse (3) reste une question ouverte.\\

\noindent \hspace{4 cm } {\it {\bf Références bibliographiques}}\\

\noindent $[1]$  A. Akcoglu. A pointwise ergodic theorem in $L_p$-$\mathrm{spaces}$.  Canad. J. 
\noindent Math. Soc. 27 (1975), 1075-1082.\\

\noindent$[2]$  A. Brunel. Théorème ergodique pour les opérateurs positifs à moyennes bornées sur
les espaces $L_p$ $(1<p<\infty)$.
Ergodic Theory Dynamical Systems. 12 (1992), 295-207.\\

\noindent$[3] $ A. Brunel. Le théorème ergodique pour les opérateurs positifs sur les espaces $L_p$
$(1<p<\infty)$ revisité. C.R. Acad. Sci. Paris, Série I 334 (2002), 205-207.\\

\noindent$[4]$  A.  Brunel, R. Emilion. Sur les opérateurs à moyennes bornées. C. R. Acad. Sci. Paris, Série I 298 (6) (1984). \\ 

\noindent$[5]$  R. Emilion. Mean bounded operators and mean ergodic theorems. J. Funct. Anal. 61 (1985), 1-14.\\

\noindent$[6]$  U. Krengel. Ergodic Theorems. De Gruyter Stud. Math., Vol. 6, 1985. \\

\noindent$[7]$  R. Sato.  On Brunel"s proof of a dominated ergodic theorem for positive linear operators on  $L_p$ $(1< p<\infty)$ (cf. $[3]$).

\vspace{0,5 cm}

\noindent Jean-Claude Lootgieter\\
Université Pierre et Marie Curie\\
 e-mail: jean-claude.lootgieter@upmc.fr\\

%$Les puissances$ %$A^n$ $ne sont pas bornées en norme.$

%,Il est très facile d'observer, à partir de (18),  qu'il existe des $f\in L^p$ pour lesquelles $A^n(f)$
%ne converge en aucun point de $\N$ tout en restant 
%bornees ou non ainsi que les moyennes de Césaro des $A^n(f)$.

%où $\alpha$ est la probabilité (6)
%et l'opérateur de translation sur $L^p$ $\mathrm{pour\,\,\,toute f et  tout k }\,\,\,\,\,\,  \mathrm \in \N,$
%\begin{equation} 
 %T(f)(k)=f(k+1)
%\end{equation}  
%%%\sum_{k=0}^{\infty} \alpha_k\vert f(k+1) \vert^p\leq C \sum_{k=1}^{\infty} \alpha_k\vert f(k) \vert^p
%\end{equation} 
%Comme  $\alpha_k \sim  c\,k^{-3/2}$, $T$ est donc bien une application linéaire positive de $L^p$ dans %$L^p$\\ Considérons le barycentre

\end{document}